\documentclass[12pt,reqno]{amsart}
\usepackage{amssymb}
\usepackage{verbatim,graphicx}
\begin{document}
\theoremstyle{plain}
\newtheorem{thm}{Theorem}[section]
\newtheorem{theorem}[thm]{Theorem}
\newtheorem{lemma}[thm]{Lemma}
\newtheorem{corollary}[thm]{Corollary}
\newtheorem{corollary and definition}[thm]{Corollary and Definition}
\newtheorem{proposition}[thm]{Proposition}
\newtheorem{example}[thm]{Example}
\theoremstyle{definition}
\newtheorem{notation}[thm]{Notation}
\newtheorem{claim}[thm]{Claim}
\newtheorem{remark}[thm]{Remark}
\newtheorem{remarks}[thm]{Remarks}
\newtheorem{conjecture}[thm]{Conjecture}
\newtheorem{definition}[thm]{Definition}
\newtheorem{problem}[thm]{Problem}
\newcommand{\Diff}{{\rm Diff}}
\newcommand{\fz}{\frak{z}}
\newcommand{\zar}{{\rm zar}}
\newcommand{\an}{{\rm an}}
\newcommand{\red}{{\rm red}}
\newcommand{\codim}{{\rm codim}}
\newcommand{\rank}{{\rm rank}}
\newcommand{\Pic}{{\rm Pic}}
\newcommand{\Div}{{\rm Div}}
\newcommand{\Hom}{{\rm Hom}}
\newcommand{\im}{{\rm im}}
\newcommand{\Spec}{{\rm Spec}}
\newcommand{\sing}{{\rm sing}}
\newcommand{\reg}{{\rm reg}}
\newcommand{\Char}{{\rm char}}
\newcommand{\Tr}{{\rm Tr}}
\newcommand{\Gal}{{\rm Gal}}
\newcommand{\Min}{{\rm Min \ }}
\newcommand{\Max}{{\rm Max \ }}
\newcommand{\soplus}[1]{\stackrel{#1}{\oplus}}
\newcommand{\dlog}{{\rm dlog}\,}    
\newcommand{\limdir}[1]{{\displaystyle{\mathop{\rm
lim}_{\buildrel\longrightarrow\over{#1}}}}\,}
\newcommand{\liminv}[1]{{\displaystyle{\mathop{\rm
lim}_{\buildrel\longleftarrow\over{#1}}}}\,}
\newcommand{\boxtensor}{{\Box\kern-9.03pt\raise1.42pt\hbox{$\times$}}}
\newcommand{\sext}{\mbox{${\mathcal E}xt\,$}}
\newcommand{\shom}{\mbox{${\mathcal H}om\,$}}
\newcommand{\coker}{{\rm coker}\,}
\renewcommand{\iff}{\mbox{ $\Longleftrightarrow$ }}
\newcommand{\onto}{\mbox{$\,\>>>\hspace{-.5cm}\to\hspace{.15cm}$}}
\catcode`\@=11
\def\opn#1#2{\def#1{\mathop{\kern0pt\fam0#2}\nolimits}}
\def\bold#1{{\bf #1}}%
\def\underrightarrow{\mathpalette\underrightarrow@}
\def\underrightarrow@#1#2{\vtop{\ialign{$##$\cr
 \hfil#1#2\hfil\cr\noalign{\nointerlineskip}%
 #1{-}\mkern-6mu\cleaders\hbox{$#1\mkern-2mu{-}\mkern-2mu$}\hfill
 \mkern-6mu{\to}\cr}}}
\let\underarrow\underrightarrow
\def\underleftarrow{\mathpalette\underleftarrow@}
\def\underleftarrow@#1#2{\vtop{\ialign{$##$\cr
 \hfil#1#2\hfil\cr\noalign{\nointerlineskip}#1{\leftarrow}\mkern-6mu
 \cleaders\hbox{$#1\mkern-2mu{-}\mkern-2mu$}\hfill
 \mkern-6mu{-}\cr}}}
\let\amp@rs@nd@\relax
\newdimen\ex@
\ex@.2326ex
\newdimen\bigaw@
\newdimen\minaw@
\minaw@16.08739\ex@
\newdimen\minCDaw@
\minCDaw@2.5pc
\newif\ifCD@
\def\minCDarrowwidth#1{\minCDaw@#1}
\newenvironment{CD}{\@CD}{\@endCD}
\def\@CD{\def\A##1A##2A{\llap{$\vcenter{\hbox
 {$\scriptstyle##1$}}$}\Big\uparrow\rlap{$\vcenter{\hbox{%
$\scriptstyle##2$}}$}&&}%
\def\V##1V##2V{\llap{$\vcenter{\hbox
 {$\scriptstyle##1$}}$}\Big\downarrow\rlap{$\vcenter{\hbox{%
$\scriptstyle##2$}}$}&&}%
\def\={&\hskip.5em\mathrel
 {\vbox{\hrule width\minCDaw@\vskip3\ex@\hrule width
 \minCDaw@}}\hskip.5em&}%
\def\verteq{\Big\Vert&&}%
\def\noarr{&&}%
\def\vspace##1{\noalign{\vskip##1\relax}}\relax\let\amp@rs@nd@&\iffalse}\fi
 \CD@true\vcenter\bgroup\relax\let\\=\cr\iffalse}\fi\tabskip\z@skip\baselineskip20\ex@
 \lineskip3\ex@\lineskiplimit3\ex@\halign\bgroup
 &\hfill$\m@th##$\hfill\cr}
\def\@endCD{\cr\egroup\egroup}
\def\>#1>#2>{\amp@rs@nd@\setbox\z@\hbox{$\scriptstyle
 \;{#1}\;\;$}\setbox\@ne\hbox{$\scriptstyle\;{#2}\;\;$}\setbox\tw@
 \hbox{$#2$}\ifCD@
 \global\bigaw@\minCDaw@\else\global\bigaw@\minaw@\fi
 \ifdim\wd\z@>\bigaw@\global\bigaw@\wd\z@\fi
 \ifdim\wd\@ne>\bigaw@\global\bigaw@\wd\@ne\fi
 \ifCD@\hskip.5em\fi
 \ifdim\wd\tw@>\z@
 \mathrel{\mathop{\hbox to\bigaw@{\rightarrowfill}}\limits^{#1}_{#2}}\else
 \mathrel{\mathop{\hbox to\bigaw@{\rightarrowfill}}\limits^{#1}}\fi
 \ifCD@\hskip.5em\fi\amp@rs@nd@}
\def\<#1<#2<{\amp@rs@nd@\setbox\z@\hbox{$\scriptstyle
 \;\;{#1}\;$}\setbox\@ne\hbox{$\scriptstyle\;\;{#2}\;$}\setbox\tw@
 \hbox{$#2$}\ifCD@
 \global\bigaw@\minCDaw@\else\global\bigaw@\minaw@\fi
 \ifdim\wd\z@>\bigaw@\global\bigaw@\wd\z@\fi
 \ifdim\wd\@ne>\bigaw@\global\bigaw@\wd\@ne\fi
 \ifCD@\hskip.5em\fi
 \ifdim\wd\tw@>\z@
 \mathrel{\mathop{\hbox to\bigaw@{\leftarrowfill}}\limits^{#1}_{#2}}\else
 \mathrel{\mathop{\hbox to\bigaw@{\leftarrowfill}}\limits^{#1}}\fi
 \ifCD@\hskip.5em\fi\amp@rs@nd@}
\newenvironment{CDS}{\@CDS}{\@endCDS}
\def\@CDS{\def\A##1A##2A{\llap{$\vcenter{\hbox
 {$\scriptstyle##1$}}$}\Big\uparrow\rlap{$\vcenter{\hbox{%
$\scriptstyle##2$}}$}&}%
\def\V##1V##2V{\llap{$\vcenter{\hbox
 {$\scriptstyle##1$}}$}\Big\downarrow\rlap{$\vcenter{\hbox{%
$\scriptstyle##2$}}$}&}%
\def\={&\hskip.5em\mathrel
 {\vbox{\hrule width\minCDaw@\vskip3\ex@\hrule width
 \minCDaw@}}\hskip.5em&}
\def\verteq{\Big\Vert&}
\def\novarr{&}
\def\noharr{&&}
\def\SE##1E##2E{\slantedarrow(0,18)(4,-3){##1}{##2}&}
\def\SW##1W##2W{\slantedarrow(24,18)(-4,-3){##1}{##2}&}
\def\NE##1E##2E{\slantedarrow(0,0)(4,3){##1}{##2}&}
\def\NW##1W##2W{\slantedarrow(24,0)(-4,3){##1}{##2}&}
\def\slantedarrow(##1)(##2)##3##4{%
\thinlines\unitlength1pt\lower 6.5pt\hbox{\begin{picture}(24,18)%
\put(##1){\vector(##2){24}}%
\put(0,8){$\scriptstyle##3$}%
\put(20,8){$\scriptstyle##4$}%
\end{picture}}}
\def\vspace##1{\noalign{\vskip##1\relax}}\relax\let\amp@rs@nd@&\iffalse}\fi
 \CD@true\vcenter\bgroup\relax\let\\=\cr\iffalse}\fi\tabskip\z@skip\baselineskip20\ex@
 \lineskip3\ex@\lineskiplimit3\ex@\halign\bgroup
 &\hfill$\m@th##$\hfill\cr}
\def\@endCDS{\cr\egroup\egroup}
\newdimen\TriCDarrw@
\newif\ifTriV@
\newenvironment{TriCDV}{\@TriCDV}{\@endTriCD}
\newenvironment{TriCDA}{\@TriCDA}{\@endTriCD}
\def\@TriCDV{\TriV@true\def\TriCDpos@{6}\@TriCD}
\def\@TriCDA{\TriV@false\def\TriCDpos@{10}\@TriCD}
\def\@TriCD#1#2#3#4#5#6{%
\setbox0\hbox{$\ifTriV@#6\else#1\fi$}
\TriCDarrw@=\wd0 \advance\TriCDarrw@ 24pt
\advance\TriCDarrw@ -1em
\def\SE##1E##2E{\slantedarrow(0,18)(2,-3){##1}{##2}&}
\def\SW##1W##2W{\slantedarrow(12,18)(-2,-3){##1}{##2}&}
\def\NE##1E##2E{\slantedarrow(0,0)(2,3){##1}{##2}&}
\def\NW##1W##2W{\slantedarrow(12,0)(-2,3){##1}{##2}&}
\def\slantedarrow(##1)(##2)##3##4{\thinlines\unitlength1pt
\lower 6.5pt\hbox{\begin{picture}(12,18)%
\put(##1){\vector(##2){12}}%
\put(-4,\TriCDpos@){$\scriptstyle##3$}%
\put(12,\TriCDpos@){$\scriptstyle##4$}%
\end{picture}}}
\def\={\mathrel {\vbox{\hrule
   width\TriCDarrw@\vskip3\ex@\hrule width
   \TriCDarrw@}}}
\def\>##1>>{\setbox\z@\hbox{$\scriptstyle
 \;{##1}\;\;$}\global\bigaw@\TriCDarrw@
 \ifdim\wd\z@>\bigaw@\global\bigaw@\wd\z@\fi
 \hskip.5em
 \mathrel{\mathop{\hbox to \TriCDarrw@
{\rightarrowfill}}\limits^{##1}}
 \hskip.5em}
\def\<##1<<{\setbox\z@\hbox{$\scriptstyle
 \;{##1}\;\;$}\global\bigaw@\TriCDarrw@
 \ifdim\wd\z@>\bigaw@\global\bigaw@\wd\z@\fi
 \mathrel{\mathop{\hbox to\bigaw@{\leftarrowfill}}\limits^{##1}}
 }
 \CD@true\vcenter\bgroup\relax\let\\=\cr\iffalse}\fi
 \tabskip\z@skip\baselineskip20\ex@
 \lineskip3\ex@\lineskiplimit3\ex@
 \ifTriV@
 \halign\bgroup
 &\hfill$\m@th##$\hfill\cr
#1&\multispan3\hfill$#2$\hfill&#3\\
&#4&#5\\
&&#6\cr\egroup%
\else
 \halign\bgroup
 &\hfill$\m@th##$\hfill\cr
&&#1\\%
&#2&#3\\
#4&\multispan3\hfill$#5$\hfill&#6\cr\egroup
\fi}
\def\@endTriCD{\egroup}
\newcommand{\sA}{{\mathcal A}}
\newcommand{\sB}{{\mathcal B}}
\newcommand{\sC}{{\mathcal C}}
\newcommand{\sD}{{\mathcal D}}
\newcommand{\sE}{{\mathcal E}}
\newcommand{\sF}{{\mathcal F}}
\newcommand{\sG}{{\mathcal G}}
\newcommand{\sH}{{\mathcal H}}
\newcommand{\sI}{{\mathcal I}}
\newcommand{\sJ}{{\mathcal J}}
\newcommand{\sK}{{\mathcal K}}
\newcommand{\sL}{{\mathcal L}}
\newcommand{\sM}{{\mathcal M}}
\newcommand{\sN}{{\mathcal N}}
\newcommand{\sO}{{\mathcal O}}
\newcommand{\sP}{{\mathcal P}}
\newcommand{\sQ}{{\mathcal Q}}
\newcommand{\sR}{{\mathcal R}}
\newcommand{\sS}{{\mathcal S}}
\newcommand{\sT}{{\mathcal T}}
\newcommand{\sU}{{\mathcal U}}
\newcommand{\sV}{{\mathcal V}}
\newcommand{\sW}{{\mathcal W}}
\newcommand{\sX}{{\mathcal X}}
\newcommand{\sY}{{\mathcal Y}}
\newcommand{\sZ}{{\mathcal Z}}
\newcommand{\A}{{\mathbb A}}
\newcommand{\B}{{\mathbb B}}
\newcommand{\C}{{\mathbb C}}
\newcommand{\D}{{\mathbb D}}
\newcommand{\E}{{\mathbb E}}
\newcommand{\F}{{\mathbb F}}
\newcommand{\G}{{\mathbb G}}
\newcommand{\HH}{{\mathbb H}}
\newcommand{\I}{{\mathbb I}}
\newcommand{\J}{{\mathbb J}}
\newcommand{\M}{{\mathbb M}}
\newcommand{\N}{{\mathbb N}}
\renewcommand{\P}{{\mathbb P}}
\newcommand{\Q}{{\mathbb Q}}
\newcommand{\R}{{\mathbb R}}
\newcommand{\T}{{\mathbb T}}
\newcommand{\U}{{\mathbb U}}
\newcommand{\V}{{\mathbb V}}
\newcommand{\W}{{\mathbb W}}
\newcommand{\X}{{\mathbb X}}
\newcommand{\Y}{{\mathbb Y}}
\newcommand{\Z}{{\mathbb Z}}
\title{A Note on an Asymptotically Good Tame Tower}   
\author{Siman Yang}
\thanks{The author was partially supported by the program for Chang Jiang Scholars and Innovative Research Team in University}
\begin{abstract}
The explicit construction of function fields tower with many rational points relative to the genus in the tower play a key role for the construction of asymptotically good algebraic-geometric codes. In 1997 Garcia, Stichtenoth and Thomas [6] exhibited two recursive asymptotically good Kummer towers over any non-prime field. Wulftange determined the limit of one tower in his PhD thesis [13]. In this paper we determine the limit of another tower [14].
\end{abstract}

\maketitle



\section{Introduction}

Let $K=\F_q$ be the finite field of cardinality $q$, and let $\mathcal{F}=(F_i)_{i\geq 0}$ be a sequence of
algebraic function fields each defined over $K$. If $F_i\varsubsetneqq F_{i+1}$ and $K$ is the full constant field for all $i\geq 0$, and  $g(F_j)>1$ for some $j\geq 0$, we call $\mathcal{F}$ a tower.

Denoted by $g(F)$ the genus of the function field $F/\F_q$ and $N(F)$ the number of $\F_q$-rational places of $F$.
It is well-known that for given genus $g$ and finite field $\F_q$, the number of $\F_q$-rational places of a function field is upper bounded due to the Weil's theorem (cf. [11]). Let $N_q(g):=\max \{N(F)|F\,\,\mbox{is a function field of genus}\, g\,\mbox{over}\,\F_q\}$ and let
\begin{equation*}
A(q)=\displaystyle\limsup_{g\rightarrow \infty}N_q(g)/g,
\end{equation*}
the Drinfeld-Vladut bound [2] provides a general upper bound of $A(q)$
\begin{equation*}
A(q)\leq \sqrt{q}-1.
\end{equation*}

Ihara[7], and Tsfasman, Vladut and Zink [12] independently showed that this bound is met when $q$ is a square by the theory of Shimura modular curves and elliptic modular curves, respectively. For non-square $q$ the exact value of $A(q)$ is unknown. Serre[10] first showed that $A(q)$ is positive for any prime power $q$
$$
A(q)\geq c\cdot \log q
$$
with some constant $c>0$ irrelevant to $q$. It was proved in [6]
that for any tower $\mathcal{F}=(F_i)_{i\geq 0}$ defined over
$\F_q$ the sequence $N(F_n)/g(F_n)_{n\geq 0}$ is convergent. We
define the limit of the tower as
$$
\lambda(\mathcal{F})=\displaystyle\lim _{i\rightarrow \infty} N(F_i)/g(F_i).
$$
Clearly, $0\leq \lambda(\mathcal{F})\leq A(q)$. We call a tower $\mathcal{F}$ asymptotically good if
$\lambda(\mathcal{F})>0$.
To be useful towards the aim of yielding asymptotically good codes, a tower must be asymptotically good.
Practical implementation of the codes also requires explicit equations for each extension step in the tower. In 1995, Garcia and Stichtenoth [4] exhibited the first explicit tower of Artin-Schreier extensions over any finite field of square cardinality which met the upper bound of Drinfeld and Vladut.
In 1997 Garcia, Stichtenoth and Thomas [6] exhibited two explicit asymptotically good Kummer towers over any non-prime field which were later generalized by Deolalikar [1].
For other explicit tame towers, readers may look at [3], [5], [9].
The two asymptotically good Kummer towers in [6] are given as below.

Let $q=p^e$ with $e>1$, and let $F_n=\F_q(x_0, \cdots, x_n)$ with
\begin{equation}
x_{i+1}^{\frac{q-1}{p-1}}+(x_i+1)^{\frac{q-1}{p-1}}=1 \,\,\,\,(i=0, \cdots, n-1).
\end{equation}
Then $\mathcal{F}=(F_0, F_1, \cdots)$ is an asymptotically  good
tower over $\F_q$ with $\lambda (\mathcal{F})\geq 2/(q-2)$.

Let  $q$ be a prime power larger than two, and let $F_n=\F_q(x_0, \cdots, x_n)$ with
\begin{equation}
x_{i+1}^{q-1}+(x_i+1)^{q-1}=1 \,\,\,\,(i=0, \cdots, n-1).
\end{equation}
Then  $\mathcal{F}=(F_0, F_1, \cdots)$ is an asymptotically good tower over $\F_{q^2}$ with $\lambda (\mathcal{F})\geq 2/(q-2)$.

Wulftange showed in [13] that $\lambda (\mathcal{F})=2/(q-2)$ for the first tower, we will show in the next section
that the limit of the second tower is also $2/(q-2)$.

\section{The limit of the tower}

\begin{lemma}
Let $F_1=K(x,y)$ defined by Eq. (2).\\
Over $K(x)$ exactly the zeroes of $x-\alpha$, $\alpha \in \F_q \backslash \{-1\}$ are ramified in $F_1$, each of ramification index $q-1$.\\
Over $K(y)$ exactly the zeroes of $y-\alpha$, $\alpha \in \F_q^*$ are ramified in $F_1$, each of ramification index $q-1$.
\end{lemma}

\begin{proof} By applying the theory of Kummer extension (cf. [11, Chap. III.7.3]).
\end{proof}
\begin{proposition}
Let $P_\alpha \in \sP (F_0)$ be a zero of $x_0-\alpha$, $\alpha \in \F_q\backslash \{-1\}$. Then, $P_\alpha$ is totally ramified in $F_{n+1}/F_n$ for any $n\geq 0$.
\end{proposition}
\begin{proof}  Let $P\in \sP (F_n)$ lying above $P_\alpha$ for some $\alpha \in \F_q\backslash \{-1\}$. From Eq. (2), one can check $x_1(P)=x_2(P)=\cdots=x_n(P)=0$. Thus the ramification index of the extension of the restriction $P$ in $K(x_i, x_{i+1})/K(x_i)$ is $q-1$ for $i=0, 1, \cdots, n$, also the ramification index of the extension of the restriction $P$ in $K(x_i, x_{i+1})/K(x_{i+1})$ is $1$ for $i=0, 1, \cdots, n$. The proof is finished by diagram chasing and repeated application of Abhyankar's lemma.
\end{proof}

Let $Q \in \sP (F_n)$ be a place ramified in $F_{n+1}$. Then $P:=Q\cap K(x_n)$ is ramified in $K(x_n, x_{n+1})$ due to Abhyankar's lemma. From Lemma 2.1, $x_n(P)=\alpha$ for some $\alpha \in \F_q\backslash \{-1\}$. If $\alpha \not= 0$, $P$ is ramifed in $K(x_{n-1}, x_n)$ of ramification index $q-1$ due to Lemma 2.1, and due to Abhyankar's lemma, the place in $K(x_{n-1},x_n)$ lying above $P$ is unramified in $K(x_{n-1}, x_n, x_{n+1})$, again by Abhyankar's lemma, $Q$ is unramified in $F_{n+1}$. Thus $Q$ is a zero of $x_n$. This implies $Q$ is a zero of $x_{n-1}-\beta$ for some $\beta \in \F_q\backslash \{-1\}$. From Eq. (2), one has the following possibilities for a place $Q \in \sP (F_n)$ ramified in $F_{n+1}$.

(a) The place $Q$ is a common zero of $x_0, x_1, \cdots, x_n$.

(b) There is some $t$, $-1\leq t <n-1$ such that

(b1) $Q$ is a common zero of $x_{t+2}, x_{t+3}, \cdots, x_n$.

(b2) $Q$ is a zero of $x_{t+1}-\alpha $ for some $\alpha \in \F_q^*\backslash \{-1\}$.

(b3) $Q$ is a common zero of $x_0 +1, x_1 +1, \cdots, x_t +1$.

 (Note that condition (b2) implies (b1) and (b3)).

\begin{lemma}
Let $-1\leq t <n$ and $Q\in \sP (F_n)$ be a place which is a zero of $x_{t+1}-\alpha$ for some $\alpha \in \F_q^*\backslash \{-1\}$. Then one has

(i)\,\, If $n<2t+2$, then $Q$ is unramified in $F_{n+1}$.

(ii) If $n\geq 2t+2$, then $Q$ is ramified in $F_{n+1}$ of ramification index $q-1$.
\end{lemma}

\begin{proof}
The assertion in (i) and (ii) follow by diagram chasing with the help of Lemma 2.1 and repeated applications of Abhyankar's lemma.
\end{proof}

For $0 \leq t <\lfloor n/2 \rfloor $ and $\alpha \in \F_q^*\backslash \{-1\}$, set

$X_{t, \alpha}:=\{ Q\in \sP (F_n)| Q \,\,\text{is a zero of}\,\, x_{t+1}-\alpha \}$
and  $A_{t, \alpha}:=\displaystyle\sum _{Q\in X_{t, \alpha}}Q$. Denote by $Q_{t+1}$ the restriction of $Q$ to $K(x_{t+1})$, we have $[F_n:K(x_{t+1})]=(q-1)^n$ and $e(Q|Q_{t+1})=(q-1)^{n-t-1}$. Then deg $A_{t, \alpha}=(q-1)^{t+1}$ follows from the fundemental equality $\sum e_i f_i =n$. Combining the above results one obtains
\begin{align}
\text{deg Diff}(F_{n+1}/F_n)&=(q-1)(q-2)+\displaystyle\sum _{\alpha \in \F_q^*\backslash \{-1\}} \displaystyle\sum _{t=0}^{\lfloor n/2 \rfloor -1}(q-2)(q-1)^{t+1}\\
&=(q-2)(q-1)^{\lfloor n/2 \rfloor +1}.
\end{align}


Now we can easily determine the genus of $F_n$ by applying the transitivity of different exponents and Hurwitz genus formula. The result is:
\[g(F_{n+1})=\left\{
\begin{array}{cccccc}
(q-2)(q-1)^{n+1}/2-(q-1)^{n/2+1}+1,\,\, \mbox{if}\,\, n\,\, \mbox{is even}, \\
(q-2)(q-1)^{n+1}/2-q(q-1)^{(n+1)/2}/2+1 ,\,\, \mbox{if}\,\, n\,\, \mbox{is odd}.
\end{array}\right.
\]

Thus $\gamma (\mathcal{F}):=\displaystyle\lim_{n\rightarrow \infty} g(F_n)/[F_n:F_0]=(q-2)/2$.

\begin{remark}
Note that from the proof of [6, Theorem 2.1 and Example 2.4],
$\gamma(\mathcal{F})$ is upper bounded by $(q-2)/2$.
\end{remark}

Next we consider the rational places in each function field $F_n$.
First we consider places over $P_{\infty}$. It is easy to see that
$P_{\infty}$ splits completely in the tower. From Prop. 2.2,
there's a unique $\F_q$-rational place in $F_n$ over $P_{\alpha}$
for any $\alpha \in \F_q\backslash \{-1\}$. Then we consider the
$K$-rational place over $P_{-1}$ in $F_n$. Let $0 \leq t <n$ and
$Q\in \sP (F_n)$ be a place which is a zero of $x_{t+1}-\alpha$
for some $\alpha \in \F_q^*\backslash \{-1\}$. We study the
condition for such place $Q$ to be $K$-rational.

\begin{lemma}
Let $Q'$ be a place of $F_2$ and $Q'$ is a zero of $x_1-\beta$ for some $\beta \in \F_q^*\backslash \{-1\}$. Then, if char$(\mathcal{F})\not=2$, $Q'$ is not a $\F_q$-rational place and  $Q'$ is a $\F_{q^2}$-rational place if and only if $\beta =-1/2$; if
char$(\mathcal{F})=2$,  $Q'$ is not a $\F_{q^2}$-rational place.

\end{lemma}
\begin{proof}
 Note that $x_2$ and $x_0+1$ both are $Q'$-prime elements. Eq. (2) implies $(\frac{x_2}{x_0+1})^{q-1}=\frac{x_1}{1+x_1}$, which is equivalent to $\beta /(1+\beta)$. $(\frac{x_2}{x_0+1})^{q-1}(Q')\not=1$ implies $Q'$ is not $\F_q$-rational, and $(\frac{x_2}{x_0+1})^{q^2-1}(Q')=1$ if and only if $\beta =-1/2$ as $\beta \in \F_q^*\backslash \{-1\}$.
\end{proof}

We generalize this result to the following proposition.

\begin{proposition}
Assume char$(\mathcal{F})$ is odd. Fix positive integers $t\leq m$. There are $2^{t-1}(q-1)$ many $\F_{q^2}\backslash \F_q$-rational places $Q$ in $\F_{q^2}(x_{m-t}, x_{m-t+1}, \cdots, x_{m+t})$ which are zeroes of $x_m-\beta$ for some $\beta \in \F_q^*\backslash \{-1\}$ if $q\equiv -1$ (mod $2^t$), with each of them corresponds to a tuple $(\alpha _1, \alpha _2, \cdots, \alpha _t)$ satisfying
\[\left\{
\begin{array}{cccccc}
x_m\equiv -1/2,\\
x_{m+1}/(x_{m-1}+1)&\equiv& \alpha _1,\,\,\mbox{with}\,\, \alpha _1 ^2&=&-1,\\
x_{m+2}/(x_{m-2}+1)&\equiv& \alpha _2,\,\,\mbox{with}\,\,\alpha _2 ^2&=&-1/\alpha_1,\\
\cdots &\equiv&\cdots, \,\, \cdots&=&\cdots,\\
x_{m+t-1}/(x_{m-t+1}+1)&\equiv& \alpha _{t-1},\,\,\mbox{with}\,\,\alpha_{t-1}^2&=&-1/\alpha_{t-2},\\
x_{m+t}/(x_{m-t}+1)&\equiv& \alpha _t,\,\,\mbox{with}\,\,\alpha _t ^{q-1}&=&-\alpha_{t-1}.
\end{array}\right.
\]
\end{proposition}

\begin{proof}
Prove by induction on $t$. For $t=1$, this is the case in Lemma 2.5, here we take $\alpha _0=1$. For $t\geq 1$, it is easily checked $(\frac{x_{m+t+1}}{x_{m-t-1}+1})^{q-1} \equiv \frac{-x_{m+t}}{x_{m-t}+1}$ from definition. Thus, $\alpha _{t+1}\in \F_{q^2}$ if and only if $\alpha_t^{q+1}=1$. By induction hypothesis on $t$, $\alpha _t ^{q-1}=-\alpha _{t-1}$. Therefore $Q$ is a $\F_{q^2}$-rational place implies $\alpha _t^2=-1/\alpha_{t-1}$. Note $\alpha_{t-1}^{2^{t-1}}=-1$. Let $q=2^tk-1$, we have $(-1)^k=1$, thus $k$ is even, i.e., $q\equiv -1$ (mod $2^{t+1}$). This finishes the induction on $t+1$.
\end{proof}

Using this proposition and Lemma 2.3, we yield the following result.
\begin{proposition}
Assume char$(\mathcal{F})$ is odd. Suppose $2^l||(q+1)$. The number of $\F_{q^2}$-rational place in $F_n$ which is a zero of  $x_m-\alpha (0<m\leq n)$ for any $\alpha \in \F_q^*\backslash \{-1\}$ is counted as below.
\[\left\{
\begin{array}{cccccc}
2^{m-1}(q-1),\,\, &\mbox{when}&\, 1\leq m \leq n/2\,\, \mbox{and}\,\, m\leq l, \\
0,\,\, &\mbox{when}&\, 1\leq m \leq n/2\,\, \mbox{and}\,\, m>l, \\
2^{n-m-1}(q-1),\,\, &\mbox{when}&\, n>m>n/2\,\, \mbox{and}\,\, n-m\leq l, \\
0,\,\, &\mbox{when}&\, n>m>n/2\,\, \mbox{and}\,\, n-m>l, \\
q-2,\,\, &\mbox{when}&\, m=n.
\end{array}\right.
\]
\end{proposition}

\begin{proof}
Let $0<m<n$ and $a=\min \{m, n-m \}$. If $a=m$ (resp. $a=n-m$), from Lemma 2.5, there exists $\F_{q^2}$-rational place in $F_{2m}$ (resp. $K(x_{2m-n}, \cdots, x_n)$) with $x_m\equiv \alpha$ for some $\alpha \in \F_q^*\backslash \{-1\}$ if and only if $2^a||(q+1)$, the number of such places is $2^{a-1}(q-1)$, and all these places totally ramified in $F_n$ according to Lemma 2.1.
\end{proof}

Hence, the number of $\F_{q^2}$-rational place in $F_n$ lying
above $P_{-1}$ is
\[\left\{
\begin{array}{cccccc}
(q-1)(2^{l+1}-1),\,\, &\mbox{if}& \,\, n>2l, \\
(q-1)(2^{(n+1)/2}-1), \,\, &\mbox{if n is odd}&\mbox{and}\, n\leq 2l, \\
(q-1)(3\times 2^{n/2-1}-1), \,\, &\mbox{if n is even}&\mbox{and}\, n\leq 2l.
\end{array}\right.
\]
\begin{remark}
If char$(\mathcal{F})>2$, among all $\F_{q^2}$-rational place in $F_n$ lying above $P_{-1}$, exactly $q-1$ are $\F_q$-rational, corresponding to  $x_n\equiv \alpha$ for some $\alpha \in \F_q^*$, respectively.
If char$(\mathcal{F})=2$, from Lemma 2.5, there are exactly $q-1$ $\F_{q^2}$-rational places in $F_n$ lying above $P_{-1}$, which are all $\F_q$-rational, corresponding to $x_n\equiv \alpha$ for some $\alpha \in \F_q^*$, respectively.
\end{remark}
Next we determine the $\F_{q^2}$-rational place $Q$ in $F_n$ lying
above $P_\alpha$ for some $\alpha \in \F_{q^2}\backslash \F_q$.
Direct calculation gives $x_1(Q)=\alpha _1$ for some $\alpha _1
\not\in \F_q$. Similarly, $x_2(Q)=\alpha -2$, $\cdots$,
$x_n(Q)=\alpha _n$, with $\alpha _i \in \overline{\F}_q \backslash
\F_q$. We observe that $Q$ is  $\F_{q^2}$-rational in $F_n$ if and
only if $\alpha, \alpha _1, \cdots, \alpha _n$ are all in
$\F_{q^2}$. To verify it, assume $\alpha, \alpha _1, \cdots,
\alpha _n$ are all in $\F_{q^2}$. Then $Q$ is completely splitting
in each extension $F_i/F_{i-1} (i=1, 2, \cdots , n)$, with
$x_i\equiv c\alpha _i$ for some $c\in \F_q ^*$ in each place
respectively.

We have $\alpha _1\in \F_{q^2}$ if and only if $(1+\alpha)^{q-1}+(1+\alpha)^{1-q}=1$. Similarly, $\alpha _i (i=1, 2, \cdots, n-1) \in \F_{q^2}$ if and only if $(1+\alpha _{i-1})^{q-1}+(1+\alpha _{i-1})^{1-q}=1$. Thus, $Q$ is $F_{q^2}$-rational implies $(1+\alpha)^{q-1}, (1+\alpha _1)^{q-1}, \cdots, (1+\alpha _{n-1})^{q-1}$ all are the root of $x^2-x+1=0$.

{\bf Claim.} $(1+\alpha)^{q-1}, (1+\alpha _1)^{q-1}, \cdots, (1+\alpha _{n-1})^{q-1}$ are equal.

{\bf Proof of the claim.} Prove by contradiction. For simplicity assume $(1+\alpha)^{q-1}\not=(1+\alpha _1)^{q-1}$.
Thus $x^2-x+1=(x-(1+\alpha)^{q-1})(x-(1+\alpha _1)^{q-1})$. Comparing the coefficient of $x^1$,
one has $1=(1+\alpha)^{q-1}+(1+\alpha _1)^{q-1}=(2+\alpha _1-\alpha _1^{q-1})/(1+\alpha_1)$. This implies $\alpha _1 \in \F_q$, which is a contradiction.

Let $p=$char$(\F_q)$, we consider the following two cases respectively.

{\bf Case 1: $p=3$.}

Since the unique root of $x^2-x+1=0$ is $-1$, $1-\alpha _1^{q-1}=-1$, and $(1+\alpha_1)^{q-1}=-1$. It is easily checked these two equalities lead to a contradiction.

{\bf Case 2: $p\not=3$.}

Thus, $-1$ is not a root of $x^2-x+1=0$, which implies $(1+\alpha)^{q-1}$ and $(1+\alpha)^{1-q}$ are distinct roots of $x^2-x+1=0$. Thus, $(1+\alpha)^{q-1}+(1+\alpha)^{1-q}=1$. By assuming $(1+\alpha)^{q-1}+(1+\alpha)^{1-q}=1$, we have $x_1(Q)=\alpha_1$, with $\alpha_1 ^{q-1}=(1+\alpha)^{1-q}$. Hence, $\alpha _1=c/(1+\alpha)$ for some $c\in \F_q^*$. Since $(1+\alpha)^{q-1}=(1+\alpha _1)^{q-1}$, direct calculation gives $c=\frac{(1+\alpha)^{2q-1}-(1+\alpha)^q}{1-(1+\alpha)^{2q-2}}$. Iterating this procedure, we have a $\F_{q^2}$-rational place $Q$ in $F_n$ lying above $P_\alpha$ for some $\alpha \in \F_{q^2}\backslash \F_q$ is one-one corresponding to a tuple $(c_1, c_2, \cdots, c_n)$ satisfying
\[\left\{
\begin{array}{cccccc}
x_0\equiv \alpha,\,\,&\mbox{with}&\,\, (1+\alpha)^{q-1}+(1+\alpha)^{1-q}=1, \\
x_1\equiv c_1/(1+\alpha):=\alpha_1 \,\, &\mbox{with}& \,\,c_1=\frac{(1+\alpha)^{2q-1}-(1+\alpha)^q}{1-(1+\alpha)^{2q-2}}\in \F_q^*, \\
x_2\equiv c_2/(1+\alpha_1):=\alpha_2 \,\, &\mbox{with}& \,\,c_2=\frac{(1+\alpha_1)^{2q-1}-(1+\alpha_1)^q}{1-(1+\alpha_1)^{2q-2}}\in \F_q^*, \\
\cdots &\mbox{with}& \,\,\cdots \\
x_{n-1}\equiv c_{n-1}/(1+\alpha_{n-2}):=\alpha_{n-1} \,\, &\mbox{with}& \,\,c_{n-1}=\frac{(1+\alpha_{n-2})^{2q-1}-(1+\alpha_{n-2})^q}{1-(1+\alpha_{n-2})^{2q-2}}\in \F_q^*, \\
x_n\equiv c_n/(1+\alpha_{n-1}),\,\, &\mbox{with}& \,\, c_n\in\F_q^*.
\end{array}\right.
\]

Therefore, for any $\alpha \in \F_{q^2}\backslash \F_q$, the number of $\F_{q^2}$-rational places in $F_n$ lying above $P_\alpha$ is zero if char$(\mathcal{F})=3$; and $(q-1)\#\{\alpha \in \F_{q^2}\backslash \F_q: (1+\alpha)^{q-1}+(1+\alpha)^{1-q}=1\}$ if char$(\mathcal{F})\not=3$.

As we have determined all $\F_{q}$-rational places and $\F_{q^2}$-rational places in $F_n$, we are now able to determine the value of $\nu(\mathcal{F})$. If char$(\mathcal{F})\not=2$ and the constant field is $\F_q$, then $\nu(\mathcal{F})=0$, and $\nu(\mathcal{F})=1$ if the constant field is $\F_{q^2}$.
If char$(\mathcal{F})=2$, and the constant field is $\F_q$ ($q>2$), then $\nu(\mathcal{F})=1$.


\begin{remark}
One can check that the function field tower recursively defined by Eq. (2) is isomorphic in some extension field of $\F_{q^2}$, to a tower recursively defined by $y^{q-1}=1-(x+\alpha)^{q-1}$, where $\alpha$ is any nonzero element of $\F_q$.
\end{remark}

From above discussion, Eq.(2) defines an asymptotically bad tower over any prime field (it does not define a tower over $\F_2$). Lenstra showed in [8, Theorem 2] that there does not exist a tower of function fields $\mathcal{F}=(F_0, F_1, \cdots)$ over a prime field which is recursively defined by $y^m=f(x)$, where $f(x)$ is a polynomial $f(x)$, $m$ and $q$ are coprime, such that the infinity place of $F_0$ splits completely in the tower, and the set $V(\mathcal{F})=\{P\in \sP(F_0)|P \,\,\mbox{is ramified in}\,\,F_n/F_0\,\, \mbox{for some}\,\, n\geq 1\}$ is finite. A tower recursively defined by Eq. (2) falls in this form with a finite set  $V(\mathcal{F})$, but no place of $F_0$ splits completely in the tower. Thus arises a problem: can one find an asymptotically good, recursive tower of the above form, over a prime field, with a finite set  $V(\mathcal{F})$ and a finite place splitting completely in the tower?

\vskip .2in
\noindent
Siman Yang\\
Department of Mathematics, East China Normal University,\\
500, Dongchuan Rd., Shanghai, P.R.China 200241. \ \ e-mail: smyang@math.ecnu.edu.cn
 \\ \\

\end{document}